\documentclass[12pt]{article}

\usepackage[latin2]{inputenc}
\usepackage{amsmath}
\usepackage{amssymb}
\usepackage{amsfonts}
\usepackage{stmaryrd}
\usepackage{theorem}
\usepackage{graphicx}
\usepackage{epstopdf}
\usepackage[parfill]{parskip} 
\usepackage{hyperref}
\hypersetup{pdftex,colorlinks=true,allcolors=blue}
\usepackage{hypcap}
\usepackage{algorithm}
\usepackage{algpseudocode} 
\usepackage{authblk}
\usepackage{tabto}

\newtheorem{thm}{Theorem}[section]
\newtheorem{prp}{Proposition}[section]
\newtheorem{lem}{Lemma}[section]

\newtheorem{dfn}{Definition}[section]
\newtheorem{ntn}{Notation}[section]

\newtheorem{prb}{Problem}[section]

\newcommand{\keywords}[1]{\par\noindent{\small{\em Keywords\/}:\ \ #1}}
\newcommand{\mscclass}[1]{\par\noindent{\small{\em MSC class\/}:\ \ #1}}

\newcommand*{\va}{\ensuremath{\alpha}}

\newcommand*{\N}{\ensuremath{\mathbb{N}}}
\newcommand*{\Z}{\ensuremath{\mathbb{Z}}}

\newcommand{\libeq}{\mathrel{\mathop:}=} 
\newcommand*{\prf}{\textbf{Proof}\ \ }
\newcommand*{\sln}{\textbf{Solution}\ \ }

\newcommand*{\sqr}{\ensuremath{\square}}
\newcommand*{\mmod}{\mathrm{mod}}

\newcommand*{\rar}{\ensuremath{\Rightarrow}}

\newcommand*{\lrar}{\ensuremath{\Leftrightarrow}}

\newcommand*{\im}{\ensuremath{\in_m}}
\newcommand*{\Zm}{\ensuremath{\mathbb{Z}_m}}
\newcommand*{\Rm}{\ensuremath{\mathrm{R}_m}}
\newcommand*{\Rme}{\ensuremath{\mathrm{R}_m^e}}
\newcommand*{\Em}{\ensuremath{\mathrm{E}_m}}
\newcommand*{\phm}{\ensuremath{\varphi(m)}}
\newcommand*{\am}{\ensuremath{|a|_m}}
\newcommand*{\bm}{\ensuremath{|b|_m}}
\newcommand*{\cm}{\ensuremath{|c|_m}}
\newcommand*{\ibma}{\ensuremath{\mathrm{ind}_b^m a}}

\newcommand*{\eibma}{\ensuremath{\exists\mathrm{ind}_b^m a}}

\newcommand*{\orbm}[1]{\ensuremath{\langle #1 \rangle_m}}
\newcommand*{\oma}{\ensuremath{\omega_m (a)}}
\newcommand*{\ind}{\ensuremath{\mathrm{ind}}}

\newcommand*{\piai}{\ensuremath{p_i^{\alpha_i}}}
\newcommand*{\Db}{\ensuremath{D_m(b,c)}}
\newcommand*{\Dc}{\ensuremath{D_m(c,b)}}

\begin{document}

\title{A Generalization of Euler's Criterion\\ to Composite Moduli}
\author{J\'{o}zsef Vass}
\affil{Dept. of Algebra and Number Theory, Institute of Mathematics\\ E\"{o}tv\"{o}s Lor\'{a}nd University\\ \vspace{0.35cm} vass@cs.elte.hu}
\date{}
\maketitle

\vspace{-0.5cm}
\begin{abstract}
\noindent A necessary and sufficient condition is provided for the solvability of a binomial congruence with a composite modulus, circumventing its prime factorization. This is a generalization of Euler's Criterion through that of Euler's Theorem, and the concepts of order and primitive roots. Idempotent numbers play a central role in this effort.\footnote{These results appeared in the author's master's thesis \cite{ma00002} and were presented at the Scientific Student Conference of E\"{o}tv\"{o}s Lor\'{a}nd University on Nov. 26, 2003 as a paper titled ``Idempotent Numbers and the Solvability of $x^k\equiv a\ (\mmod\ m)$''. The paper was renamed and some revisions were made in May 2015.}
\vspace{0.35cm}
\mscclass{\hspace{0.00cm} 11A15 (primary); 11A07, 11C08 (secondary).}
\keywords{\hspace{0.15cm} binomial congruences, power residues, generalized primitive roots.}
\end{abstract}

\vspace{-0.75cm}
\tableofcontents
\newpage

\section{Introduction} \label{s01}

\subsection{Overview} \label{s0101}

The solvability of binomial congruences of the form $x^k\equiv a\ (\mmod\ m),\ k\in\N,\ a\in\Z_m$ where the modulus $m$ is any integer, is generally reduced using the Chinese Remainder Theorem to a system of congruences with prime power moduli, for which solvability can be decided with well-known techniques. Since the algorithmic complexity of prime factorization is high, it may be worthwhile to explore an alternative path.

This path will be set by idempotent numbers $e^2\equiv e\ (\mmod\ m)$ which are projections to divisors of $m$ sharing the same prime power factors, enabling us to bypass the Chinese Remainder Theorem. Their relevance will emerge with the generalization of Euler's Theorem which becomes the basis for the concepts of order, orbit, and index. A useful generalization of primitive roots is subsequently suggested. The mentioned alternative path must somehow avoid the fact that genuine primitive roots which generate all coprime residues do not exist for a general modulus. Indeed this is accomplished with a critical theorem, leading to a theoretical equivalence condition for the solvability of such a congruence, similar to Euler's Criterion. Such criteria for power residues may lead to practical reciprocity laws.

For an overview of congruences see Andrews \cite{bb00020}, and of reciprocity see Lemmermeyer \cite{bb00021}. For a more complete discussion of composite moduli via idempotent numbers, see Vass \cite{ma00002}.

\subsection{Preliminaries} \label{s0102}

\begin{ntn} \label{s010201}
Let $\N$ denote the set of integers greater than or equal to $1$. Let the prime numbers be denoted as $p_i,\ i\in\N$ in ascending order. Denote the prime factorization of $m\in\N$ as $m=p_1^{\alpha_1}\ldots p_i^{\alpha_i}\ldots\ (\alpha_i\ge 0)$. Denote $\Zm\libeq\{1,\dots,m\}$ and let $a\ \mmod\ m$ be the number $b\in\Zm$ for which $a\equiv b\ (\mmod\ m)$. For $A\subset\Zm,\ a\in\Z$ write $a\im A$ iff $(a\ \mmod\ m) \in A$. Let $(a,b)$ denote the greatest common divisor of the numbers $a, b\in\N$. For $A\subset\N$ let $\gcd(a: a\in A)$ denote the greatest common divisor of all the elements in $A$. Let $[a,b]$ denote the least common multiple. Let $\phm$ denote Euler's totient function.
\end{ntn}

\begin{thm} \label{s010202} \textup{(Euler's Theorem \cite{ba00119})}\
$\forall m\in\N,\ a\in\Zm,\ (a,m)=1:\ a^{\phm}\equiv 1\ (\mmod\ m)$.
\end{thm}

\begin{thm} \label{s010203} \textup{(Euler's Criterion \cite{ba00117, ba00118})}\
Take a modulus $m$ of the form $2, 4, p^\va$ or $2p^\va$ with an odd prime number $p$ and $\va\in\N$ (i.e. a primitive root exists). Then $a\in\Zm,\ (a,m)=1$ is a $k$-th power residue ($k\in\N$), meaning $x^k\equiv a\ (\mmod\ m)$ is solvable for $x\in\Zm$ iff
\[ a^{\frac{\phm}{(k, \phm)}}\equiv 1\ (\mmod\ m). \]
\end{thm}

The proof of the above criterion relies heavily on the existence of a primitive root for moduli of the above form. So to find a criterion for other moduli, the challenge becomes to avoid the need for a primitive root.

\section{Idempotent and Regular Numbers} \label{s02}

\subsection{Order} \label{s0201}

\begin{dfn} \label{s020101}
A residue $e\in\Zm$ is an \textbf{idempotent number} modulo $m$ if $e^2\equiv e\ (\mmod\ m)$, and let $\Em$ denote their set.
\end{dfn}

It is easy to show that their cardinality is $|\Em|= 2^N$ where $N$ is the number of distinct prime power factors of $m$ (so if $m$ is a prime power, then $\Em=\{1,m\}$). The notation $e$ comes from the first letter of the Hungarian word for ``unit'', since as stated in Theorem \ref{s020206} certain subsets of $\Zm$ form abelian groups with an idempotent number as their unit element.

\begin{thm} \label{s020102} \textup{(Generalization of Euler's Theorem)}\
$\forall m\in\N,\ a\in\Zm:\ a^{\phm}\im \Em$.
\end{thm}
\noindent
\prf
Take any $i\in\N$ index for which $\alpha_i>0$ in the prime factorization of $m$. Let us consider two cases, depending on whether $p_i$ divides $a$ or not. Supposing first that $p_i\mid a$
\[ \alpha_i = 1+(\alpha_i -1)\le 2^{\alpha_i -1}\le p_i^{\alpha_i -1}\le p_i^{\alpha_i -1}(p_i -1) \le \phm \]
we may conclude that $a^{\phm}\equiv 0\ (\mmod\ \piai)$. On the other hand, if $p_i\nmid a$ then by Euler's Theorem \ref{s010202} and $\varphi(p_i^{\alpha_i})\mid\phm$ we get that $a^{\phm}\equiv 1\ (\mmod\ \piai)$. Thus in both cases
\[ a^{\phm}(a^{\phm}-1)\equiv 0\ (\mmod\ \piai) \]
for any $i$ index, implying that $a^{\phm}\ \mmod\ m$ is idempotent.\ \sqr

\begin{dfn} \label{s020103}
For $a\in\Z$ let its \textbf{order} modulo $m$ be the smallest $n\in\N$ power for which $a^n \im \Em$. Let $\am$ denote this $n$ which exists due to the above theorem.
\end{dfn}

\subsection{Regularity} \label{s0202}

\begin{dfn} \label{s020201}
The residue $a\in\Zm$ is said to be \textbf{regular} modulo $m$ if $a^{\am +1}\equiv a\ (\mmod\ m)$ and let $\Rm$ denote their set. For $e\in\Em$ denote $\Rme\libeq \{a\in\Rm:\ a^{\am}\equiv e\ (\mmod\ m)\}$.
\end{dfn}

Among many interesting facts, it is true that all residues are regular modulo $m$ iff $m$ is square-free. Several equivalent definitions may be given for regularity. Perhaps the most straightforward one is that $a$ is regular iff there exists some power $n>1$ for which $a^n$ is congruent to $a$. In essence, $a\in\Rm$ iff $p_i\mid a$ implies $p_i^{\alpha_i}\mid a$. Note also that $\Rm^1$ is a reduced residue system modulo $m$. (See the author's master's thesis \cite{ma00002} for the proofs.)

\begin{prp} \label{s020202}
For any $a\in\Rm,\ k,l\in\N$ the following hold:
\begin{enumerate}
\item $a^k\im\Em\ \rar\ \am\mid k$,
\item $\am\mid\phm$,
\item $a^k\equiv a^l\ (\mmod\ m)\ \lrar\ k\equiv l\ (\mmod\ \am)$,
\item $|a^k|_m = \am/(k,\am)$.
\end{enumerate}
\end{prp}
\noindent
\prf
1. Let $q,r\in\N\cup\{0\}$ be such that $k=q\am +r,\ 0\le r<\am$. Then
\[ a^k\equiv (a^{\am})^q\cdot a^r\equiv a^{\am}\cdot a^r\equiv a^r\ (\mmod\ m) \]
so $a^r\im\Em$, which can only be if $r=0$, by the definition of order.\\
2. Follows from 1.\\
3. Clearly we have
\[ a^k\equiv a^l\ (\mmod\ m)\ \rar\ a^k a^{l\phm -l}\equiv a^{l\phm}\ (\mmod\ m). \]
Since $a^{l\phm}\im\Em$ then by 1. and 2. we have
\[ 0\equiv k+l\phm -l\equiv k-l\ (\mmod\ \am)\ \rar\ k\equiv l\ (\mmod\ \am). \]
Now if $l\ge k$ and $k\equiv l\ (\mmod\ \am)$, then for some $q\ge 0$, we have $l=k+q\am$, so
\[ a^l\equiv a^{k+q\am}\equiv a^k a^{\am}\equiv a^k\ (\mmod\ m) \]
where the last congruence holds, because $a$ is regular.\\
4. Considering the congruence
\[ (a^k)^{\frac{\am}{(k,\am)}} = (a^{\am})^{\frac{k}{(k,\am)}}\im\Em \]
we have $|a^k|_m\le \am/(k,\am)$ by the definition of order. Also by 1. we have
\[ a^{kl}\im\Em\ \rar\ \am\mid kl\ \lrar\ \frac{\am}{(k,\am)}\mid l \]
so we have $|a^k|_m\ge \am/(k,\am)$.\ \sqr

\begin{prp} \label{s020204}
A number $a\in\Zm$ is regular iff the following equivalence holds
\[ a^k\equiv a^l\ (\mmod\ m)\ \lrar\ k\equiv l\ (\mmod\ \am)\ \ (k,l\in\N). \]
\end{prp}
\noindent
\prf
By Proposition \ref{s020202}, we have that if $a$ is regular, then the equivalence holds. On the other hand, if the equivalence holds, then with $k\libeq\am +1,\ l\libeq 1$ we have that $a$ is regular.\ \sqr

\begin{dfn} \label{s020205}
Denote $a^0\libeq a^{\am}\ \mmod\ m$. Let the \textbf{inverse} of $a\in\Rm$ be the residue $a^{-1}\libeq a^{\am -1}\ \mmod\ m$, and for any $n\in\N$ denote $a^{-n}\libeq (a^{-1})^n\ \mmod\ m$.
\end{dfn}

\begin{thm} \label{s020206}
For all $e\in\Em$ the structure $\langle \Rme; \{e,^{-1},\cdot\} \rangle$ is an
abelian group.
\end{thm}
\noindent
\prf
The properties to be shown are mostly trivial, except for maybe one. We need to show
that for all $a\in\Rme$ there exists a unique $b\in\Rme$ such that
$ab\equiv e\ (\mmod\ m)$.\\
Let $b\libeq a^{\am -1}\ \mmod\ m$. It is obvious that $ab\equiv e\ (\mmod\ m)$. Now, let us suppose that there exists some other $b'\in\Rme$ such that $ab'\equiv e\ (\mmod\ m)$. Then we have
\[ a(b-b')\equiv 0\ (\mmod\ m)\ \rar\ 0\equiv a^{\am -1}\cdot a(b-b')\equiv \]
\[ \equiv e(b-b')\equiv b^{\bm +1}-(b')^{|b'|_m +1}\equiv b-b'\ (\mmod\ m).\ \sqr \]

\begin{prp} \label{s020207}
For $a\in\Rm,\ n\in\N,\ i,j\in\Z$
\[ (a^n)^{-1} \equiv a^{-n}\ (\mmod\ m) \]
\[ a^{i+j} \equiv a^i\cdot a^j\ (\mmod\ m). \]
\end{prp}
\noindent
\prf
The first statement is equivalent to saying that
\[ (a^n)^{|a^n|_m-1}\equiv a^{n\am -n}\ (\mmod\ m) \]
which by Proposition \ref{s020204} is equivalent to (when $n\frac{\am}{(n,\am)}-n\ne 0$)
\[ \frac{n}{(n,\am)}\am-n\equiv n\am -n\ (\mmod\ \am) \]
and this congruence clearly holds. In the omitted case
\[ n\frac{\am}{(n,\am)}-n = 0\ \lrar\ \am\mid n \]
so for some $k\in\N$, we have
\[ a^{-n}\equiv a^{n\am -n} = a^{(n-k)\am}\equiv a^0\equiv (a^n)^{-1}\ (\mmod\ \am). \]
For the second property, we can distinguish four different cases (for nonzero exponents):\\
The case of $i,j>0$ is trivial. The case of $i,j<0$:
\[ a^{i+j}=a^{-|i+j|}\equiv (a^{-1})^{|i+j|}=(a^{-1})^{|i|}\cdot (a^{-1})^{|j|}\equiv \]
\[ \equiv a^{-|i|}\cdot a^{-|j|}\equiv a^i\cdot a^j\ (\mmod\ m). \]
The case of $j\ge|i|$:
\[ a^{i+j}=a^{j-|i|}\ \rar\ a^j=a^{i+j}\cdot a^{|i|}\ \rar \]
\[ \rar a^{i+j}\equiv a^j\cdot (a^{|i|})^{-1}\equiv a^j\cdot a^{-|i|}=a^i\cdot a^j\ (\mmod\ m). \]
The case of $j<|i|$:
\[ a^{i+j}\equiv a^{j-|i|}\equiv a^{-(|i|-j)}\equiv (a^{|i|-j})^{-1}\equiv (a^{|i|}\cdot a^{-j})^{-1}\ (\mmod\ m) \]
where the last congruence is true with the application of the previous case. Lastly\\
\[ (a^{|i|}\cdot a^{-j})\cdot (a^{-|i|}\cdot a^j)\equiv (a^{|i|}) (a^{|i|})^{-1} (a^j)^{-1} (a^j)\equiv (a^{\am})^{|i|+j}\equiv a^{\am}\ (\mmod\ m) \]
so by the unicity of the inverse (previous theorem), we have
\[ (a^{|i|}\cdot a^{-j})^{-1}\equiv a^{-|i|}\cdot a^j\equiv a^i\cdot a^j\ (\mmod\ m). \]
The case of $i>0,\ j<0$ is similar to the previous two.\ \sqr

\subsection{Orbit} \label{s0203}

\begin{dfn} \label{s020301}
Let the \textbf{orbit} of $a\in\Zm$ be the set $\orbm{a}\libeq \{a^n\ \mmod\ m:\ 1\le n\le \am\}$.
\end{dfn}

\begin{prp} \label{s020302}
For any $b,c\in\Rm,\ n,k\in\N$ we have
\[ b^n,b^k\im\orbm{c}\ \lrar\ b^{(n,k)}\im\orbm{c}. \]
\end{prp}
\noindent
\prf
First suppose that $b^n\equiv c^i,\ b^k\equiv c^j\ (\mmod\ m)$. Without hurting generality, we may suppose that there exist $x,y\ge 0$ such that $(n,k)=nx-ky$. So we have
\[ b^{(n,k)}=b^{nx-ky}=b^{nx+(-ky)}\equiv b^{nx}\cdot b^{-ky}\equiv b^{nx}\cdot (b^{ky})^{-1}\equiv (c^{ix})\cdot (c^{jy})^{\phm -1}\im\orbm{c} \]
with the application of Proposition \ref{s020207}.\\
Now, let us suppose that $b^{(n,k)}\equiv c^l\ (\mmod\ m)$. Then we have
\[ b^n\equiv b^{(n,k)\frac{n}{(n,k)}}\equiv (c^l)^{\frac{n}{(n,k)}}\im\orbm{c} \]
and also $b^k\im\orbm{c}$ similarly.\ \sqr

\begin{dfn} \label{s020303}
For $e\in\Em,\ b,c\in\Rme$, denote
\[ D_m(b,c)\libeq \gcd(n\in\N:\ 1\le n\le\bm,\ b^n\im\orbm{c}). \]
\end{dfn}

\begin{prp} \label{s020304}
If $e\in\Em,\ b,c\in\Rme$, then $D_m(b,c)\mid\bm$ and
\[ b^k\im\orbm{c}\ \ \lrar\ \ D_m(b,c)\mid k. \]
Furthermore $b^{D_m(b,c)}\im\orbm{c}$ and
\[ \orbm{b}\cap\orbm{c} = \orbm{b^{D_m(b,c)}}\ \ \mathrm{and}\ \ |\orbm{b}\cap\orbm{c}|=\frac{\bm}{D_m(b,c)}. \]
\end{prp}
\noindent
\prf
By the previous theorem and induction $b^{D_m(b,c)}\im\orbm{c}$. Supposing that $D_m(b,c)\mid k$ we have
\[ b^k\equiv (b^{D_m(b,c)})^{\frac{k}{D_m(b,c)}}\im\orbm{c}. \]
If $b^k\im\orbm{c}$ then with $k'\libeq k\ \mmod\ \bm$ we have $b^{k'}\im\orbm{c}$ so $D_m(b,c)\mid k'$ by definition, and from this it follows that $D_m(b,c)\mid k$.\\
By the first property now proven, we get the second one
\[ \orbm{b}\cap\orbm{c} = \orbm{b^{D_m(b,c)}}. \]
It is also true that $D_m(b,c)\mid\bm$ since
\[ b^{\bm}\equiv e\equiv c^{\cm}\im\orbm{c} \]
so lastly, we have that
\[ |\orbm{b}\cap\orbm{c}|=|\orbm{b^{D_m(b,c)}}|=|b^{D_m(b,c)}|_m = \frac{\bm}{(D_m(b,c),\bm)}=\frac{\bm}{D_m(b,c)}.\ \sqr \]

\subsection{Index} \label{s0204}

\begin{dfn} \label{s020401}
If it exists for $a,b\in\Zm$, let the \textbf{index} \ibma\ denote the smallest $n\in\N$, for which $b^n\equiv a\ (\mmod\ m)$. Let this existence be denoted as \eibma. For $a\in\Rm$ let its \textbf{primitive order} be the number $\oma\libeq \max\{\bm:\ b\in\Rm,\ \eibma\}$.
\end{dfn}

If $(a,m)=1$ and a primitive root exists modulo $m$, then clearly $\oma=\phm=|g|_m$ for any primitive root $g\in\Rm^1$. Thus a number $g\in\Rm$ may be considered a ``generalized primitive root'' if $\omega_m(g)=|g|_m$ (see \cite{ma00002} for further discussion).

\begin{prp} \label{s020403}
For any $k\in\N,\ e\in\Em,\ a,b\in\Rme,\ \eibma$ we have the equivalence
\[ (k,\bm)\mid\ibma\ \ \lrar\ \ a^{\frac{\bm}{(k,\bm)}}\im\Em. \]
\end{prp}
\noindent
\prf The equivalence can be deduced as follows.
\[ e\equiv a^{\frac{\bm}{(k,\bm)}}\equiv b^{\ibma\frac{\bm}{(k,\bm)}}\equiv b^{\bm\frac{\ibma}{(k,\bm)}}\ (\mmod\ m) \]
\[ \lrar\ \bm\mid \bm\frac{\ibma}{(k,\bm)}\ \lrar\ \frac{\ibma}{(k,\bm)}\in\N\ \lrar\ (k,\bm)\mid\ibma.\ \sqr \]

\begin{prp} \label{s020404}
If $e\in\Em,\ a,b\in\Rme,\ (\am,\bm)=1$ then $|ab|_m=\am\cdot\bm$.
\end{prp}
\noindent
\prf
We readily see that $(ab)^{\am\cdot\bm}\equiv e\ (\mmod\ m)$ implying $|ab|_m\mid\am\cdot\bm$. For the other direction of division, we first deduce
\[ e\equiv (ab)^{\am\cdot |ab|_m}\equiv e\cdot b^{\am\cdot |ab|_m}\equiv b^{\am\cdot |ab|_m}\ (\mmod\ m)\ \rar\ \bm\mid \am\cdot |ab|_m\ \rar\ \bm\mid |ab|_m \]
and similarly $\am\mid |ab|_m$ also holds, implying that $\am\cdot\bm\mid |ab|_m$.\ \sqr

\begin{lem} \label{s020406}
Given $u,v,w\in\N,\ w\mid (u,v)$ there exist $u_{1,2}, v_{1,2}, w_{1,2}\in\N$ such that $u= u_1 u_2,\ v= v_1 v_2,\ w= w_1 w_2$ and $(u,v)= u_2 v_1$ and $w_1\mid v_1\mid u_1,\ w_2\mid u_2\mid v_2$ and $1= (u_1,u_2)= (v_1,v_2)= (w_1,w_2)= (u_1,v_2)= (u_2,v_1)$.
\end{lem}
\noindent
\prf
Letting $C\libeq (u,v),\ U\libeq u/C,\ V\libeq v/C$ we have $(U,V)=1$. Partitioning $C$ according to the prime factors of $U$ and $V$, there must exist $A,B\in\N\ (C=AB)$ such that $(A,B)=1= (A,V)= (B,U)$. Clearly $u=AUB,\ v=AVB$ so defining $u_1\libeq AU,\ u_2\libeq B,\ v_1\libeq A,\ v_2\libeq VB$ then due to $w\mid C=AB= u_2 v_1$ there must exist $w_{1,2}\in\N\ (w= w_1 w_2)$ such that $w_1\mid v_1,\ w_2\mid u_2$ and clearly $v_1\mid u_1,\ u_2\mid v_2$. Lastly, observe that $1= (u_1,u_2)= (v_1,v_2)= (w_1,w_2)= (u_1,v_2)= (u_2,v_1)$ as required.\ \sqr

This lemma resembles Kalm\'{a}r's Four-Number Theorem \cite{ba00120} which can be employed to show the Fundamental Theorem of Arithmetic, while bypassing the need for the concepts of the ``greatest common divisor'' or the ``least common multiple'', which are two typical approaches. Similarly, our quest to show a generalization of Euler's Criterion hinges on this lemma and the theorem below to be shown with it, bypassing this time the lack of a cyclical generator (a ``genuine'' primitive root) for most composite moduli.

\begin{thm} \label{s020405}
Suppose that $e\in\Em,\ a,b,c\in\Rme$ and $a\in\orbm{b}\cap\orbm{c}$. Then there exists some $d\in\Rme$ for which $a\in\orbm{d}$ and $|d|_m=[\bm,\cm]$.
\end{thm}
\noindent
\textbf{Proof}\footnote{The theorem was conjectured by the author, and the presented proof is a slightly modified version of the one provided by Prof. Mih\'{a}ly Szalay.}
By Proposition \ref{s020304} we have
\[ \orbm{b}\cap\orbm{c}=\orbm{b^{D_m(b,c)}}=\orbm{c^{D_m(c,b)}} \]
so there exists some $K\in\N$ such that
\[ (b^{D_m(b,c)})^K\equiv c^{D_m(c,b)}\ (\mmod\ m). \]
Therefore from
\[ |b^{\Db}|_m=|\orbm{b}\cap\orbm{c}|=|c^{\Dc}|_m=\frac{|b^{\Db}|_m}{(K,|b^{\Db}|_m)} \]
we get that $(K,|b^{\Db}|_m)=1$. Furthermore
\[ \frac{\bm}{\Db}=|\orbm{b}\cap\orbm{c}|=\frac{\cm}{\Dc}\ \rar\ \Dc\frac{\bm}{(\bm,\cm)}=\Db\frac{\cm}{(\bm,\cm)} \]
\[ \rar\ \frac{\bm}{(\bm,\cm)}\mid\Db\frac{\cm}{(\bm,\cm)}\ \mathrm{and\ since}\ \left(\frac{\bm}{(\bm,\cm )},\ \frac{\cm}{(\bm,\cm )}\right)=1\ \rar \]
\[ \frac{\bm}{(\bm,\cm)}\mid\Db\ \mathrm{and}\ w\mid (\bm,\cm)\ \mathrm{with}\ w\libeq\frac{\Db (\bm,\cm)}{\bm}\in\N. \]
According to Lemma \ref{s020406}, for $u\libeq\bm,\ v\libeq\cm$ the following factorization is possible
\[ \bm= u_1 u_2,\ \ \cm= v_1 v_2,\ \ w= w_1 w_2 \mid (\bm,\cm)= u_2 v_1 \]
\[ w_1\mid v_1\mid u_1,\ \ w_2\mid u_2\mid v_2,\ \ 1= (u_1,u_2)= (v_1,v_2)= (u_1,v_2)= (u_2,v_1). \]
Then these properties hold
\[ |b^{u_2}|_m = \frac{\bm}{(u_2,\bm)}= u_1,\ \ |c^{v_1}|_m= \frac{\cm}{(v_1,\cm)}= v_2,\ \ (|b^{u_2}|_m, |c^{v_1}|_m)=1 \]
\[ \Db= \frac{w\bm}{(\bm,\cm)}= w \frac{u_1 u_2}{u_2 v_1}= w\frac{u_1}{v_1},\ \ \Dc=\Db\frac{\cm}{\bm}=w\frac{u_1}{v_1}\frac{v_1 v_2}{u_1 u_2}= w\frac{v_2}{u_2} \]
\[ |b^{\Db}|_m= \frac{\bm}{(\Db,\bm)}= \frac{u_1 u_2}{\Db}=\frac{u_1 u_2}{w\frac{u_1}{v_1}}= \frac{u_2 v_1}{w}\in\N. \]
Defining $d\libeq  b^{u_2} c^{v_1}\ \mmod\ m$ we have by Proposition \ref{s020404} the required order
\[ |d|_m = u_1 v_2= \frac{u_1 u_2 v_1 v_2}{u_2 v_1}= \frac{\bm\cm}{(\bm,\cm)}= [\bm,\cm]. \]
Lastly, we need an exponent $E\in\N$ such that $d^E\equiv a\ (\mmod\ m)$. Defining
\[ E\libeq w\frac{u_1 v_2}{v_1 u_2}\ \ \rar\ \ d^E\equiv (b^{\Db})^{v_2}(c^{\Dc})^{u_1}\equiv (b^{\Db})^{v_2+ K u_1}\ (\mmod\ m). \]
Now observe that $(v_2+ K u_1, |b^{\Db}|_m)=1$ where $|b^{\Db}|_m=u_2 v_1/w$ from above, since $\frac{v_1}{w_1}\mid v_1\mid u_1\mid K u_1$ but $(v_1/w_1, v_2)=1$ and $\frac{u_2}{w_2}\mid u_2\mid v_2$ but $(u_2/w_2, K u_1)=1$ since as we saw above $1=(K,|b^{\Db}|_m)= (K, u_2 v_1/w)$. So there must exist an inverse $N\in\N$ such that $(v_2+ K u_1)N\equiv 1\ (\mmod\ |b^{\Db}|_m)$. Furthermore, by the assumption of the theorem, there exists an $I\in\N$ such that $(b^{\Db})^I\equiv a\ (\mmod\ m)$.\\
Multiplying the above exponents, we may now conclude that $a\in\orbm{d}$ since
\[ d^{ENI}\equiv (b^{\Db})^{(v_2+ K u_1)NI}\equiv (b^{\Db})^I\equiv a\ (\mmod\ m).\ \sqr \]

\section{Solvability} \label{s03}

\begin{prp} \label{s030102}
For any $m\in\N,\ a\in\Rm,\ k\in\N$, if the equation $x^k\equiv a\ (\mmod\ m)$ is solvable for $x\in\Zm$ then necessarily
\[ a^{\frac{\phm}{(k, \phm)}}\im\Em. \]
\end{prp}
\noindent
\prf
Letting one of the solutions be denoted as $x_0$ we have
\[ a^{\frac{\phm}{(k,\phm)}} \equiv (x_0^k)^{\frac{\phm}{(k,\phm)}} \equiv (x_0^{\phm})^{\frac{k}{(k,\phm)}} \im\Em.\ \sqr \]

\begin{thm} \label{s030101} \textup{(Generalization of Euler's Criterion)}\
For any $m\in\N,\ a\in\Rm,\ k\in\N$ the equation $x^k\equiv a\ (\mmod\ m)$ is solvable for $x\in\Zm$ if and only if
\[ a^{\frac{\oma}{(k, \oma)}}\im\Em. \]
\end{thm}
\noindent
\prf
Let $b\in\Rm$ be such that $\eibma$ and $\bm=\oma$. Then by Proposition \ref{s020403}
\[ a^{\frac{\oma}{(k,\oma)}}\im\Em\ \lrar\ (k,\bm)\mid\ibma. \]
If $(k,\bm)\mid\ibma$ holds, then there exists some $1\le l\le\bm$ for which $kl\equiv \ibma\ (\mmod\ \bm)$. Therefore
\[ b^{kl}\equiv b^{\ibma}\ (\mmod\ m)\ \rar\ (b^l)^k\equiv a\ (\mmod\ m) \]
implying that $b^l$ is a solution of the equation. Conversely, suppose that $x_0$ is a solution, and denote $e\libeq a^{\am}\ \mmod\ m,\ c\libeq x_0 e\ \mmod\ m$. Then $c$ must be a regular solution, since
\[ c^k\equiv (x_0)^k e\equiv a\cdot a^{\am}\equiv a\ (\mmod\ m) \]
\[ c\cdot c^{\cm}\equiv c\cdot c^{\phm}\equiv x_0 e (x_0 e)^{\phm}\equiv x_0 e (x_0^{\phm})^k\equiv x_0 e (x_0^k)^{\phm}\equiv x_0 e\equiv c\ (\mmod\ m). \]
We now show that $\cm\mid\bm$. Supposing indirectly that $\cm\nmid\bm$ we have $\cm < \bm$ by the definition of $\oma=\bm$. We also know by Theorem \ref{s020405} that there exists some $d\in\Rm$ such that $\exists\ind_d^m a$ and $|d|_m=[\bm,\cm]$. Then $\cm\nmid\bm$ implies that $|d|_m>\bm$ which contradicts our original selection of $b$. So we must have that $\cm\mid\bm$ implying
\[ a^{\frac{\oma}{(k,\oma)}}\equiv a^{\frac{\bm}{(k,\bm)}}\equiv (c^k)^{\frac{\bm}{(k,\bm)}}\equiv (c^{\cm})^{\frac{\bm}{\cm}\cdot\frac{k}{(k,\bm)}} \equiv e\ (\mmod\ m).\ \sqr \]

\section{Concluding Remarks} \label{s04}

A generalization of Euler's Criterion was presented in Theorem \ref{s030101}, while the lack of a cyclical generator (primitive root) in general, was circumvented via Theorem \ref{s020405}. The criterion
\[ a^{\frac{\oma}{(k,\oma)}}\im\Em \]
in its current form is theoretical. For its practical verification, the calculation of $\oma$ must be made efficient. Likely the examination of the mapping $m\mapsto \oma$ is a worthwhile direction for future investigations, since $\oma=\phm$ when $(a,m)=1$ and a primitive root exists modulo $m$.

This paper was inspired by the following solution devised by the author, upon accidentally employing Euler's Theorem when $(a,m)\neq 1$ and seeing that $a^{\phm}\ \mmod\ m$ is idempotent. This problem can nevertheless be solved in an elementary way as well.

\begin{prb} \label{s040101}
Defining the sequence of numbers $(a_n)$ recursively as
\[ a_0\libeq 1,\ a_{n}\libeq 42^{a_{n-1}}\ (n\in\N) \]
what are the last two digits of $a_{100}$?
\end{prb}
\noindent
\sln
Let us first calculate the order and idempotent number for the last few terms, where each modulus is implied by the previous order. We descend in modulus until reaching the term $a_{97}$ congruent to zero -- this must necessarily occur since $|a|_m\leq\phm< m$.
\[ a_{100}=42^{a_{99}},\ |42|_{100}=20,\ 42^{20}\equiv 76\ (\mmod\ 100) \]
\[ a_{99}=42^{a_{98}},\ |42|_{20}=4,\ 42^4\equiv 16\ (\mmod\ 20) \]
\[ a_{98}=42^{a_{97}},\ |42|_4=2,\ 42^2\equiv 0\ (\mmod\ 4). \]
We reach zero with $a_{97}\equiv 0\ (\mmod\ 2)$ since $2\mid 42\mid a_{97}$, implying $a_{97}=2i,\ i\in\N$. Now working backwards
\[ a_{98}=42^{2i}\equiv 0\ (\mmod\ 4)\ \rar\ a_{98}=4j,\ j\in\N \]
\[ a_{99}=42^{4j}\equiv 16\ (\mmod\ 20)\ \rar\ a_{99}=20k+16,\ k\in\N \]
we finally arrive at
\[ a_{100}= 42^{20k+16}\equiv 76\cdot 42^{16}\equiv 76\cdot 56\equiv 56\ (\mmod\ 100).\ \sqr \]

The author is grateful to Prof. Mih\'{a}ly Szalay for providing the proof of Theorem \ref{s020405} and for his careful review of this paper, as well as for that of Prof. Andr\'{a}s S\'{a}rk\"{o}zy.

\bibliographystyle{abbrv}
\bibliography{mybib}
\addcontentsline{toc}{section}{\textbf{References}}

\end{document}